# Approximate solution to the fractional Lane-Emden type equations

M. I. Nouh[1,2] and Emad A-B. Abdel-Salam[3]

[1]Department of Physics, College of Science, Northern Border University, Arar, Saudi Arabia.

Email: nouh@nbu.edu.sa

[2]Department of Astronomy, National Research Institute of Astronomy and Geophysics, Helwan, Cairo17211, Egypt.

[3]Department of Mathematics, Faculty of Science, Assiut University, New Valley Branch, El-Kharja 72511, Egypt.

Email: emad.mahmoud@newvedu.au.edu.eg

**Abstract.**

In this paper, approximate solutions for a class of fractional Lane–Emden type equations based on the series expansion method are presented. Various examples are introduced and discussed. The recurrence relation for the components of the approximate solution is constructed. In the standard integer order derivative, it is shown that the results are the same as those obtained by Adomian decomposition method, homotopy perturbation method, modified Laplace decomposition method and variational iteration method.

1 **Introduction**

The Lane–Emden equation (LEE) may be considered as a Poisson equation for the gravitational potential of self-gravitating, spherically symmetric polytropic field Dehghan, Shakeri (2008).

Numerical integration is the common method used to solve LEE. In the other hand, there is many approximate methods to solve the Lane–Emden equation [1]. The first method is to transform the equation into an integrodifferential equation and then iterate this equation [2]. The second method is developed by [1], [3] and [4]. This method uses the Adomian decomposition method [5] and introduces the solution in the form of a power series. The third method is to assume a power series solution and use an accelerating techniques to obtain a convergent series [6].

Various real life problems are described exactly based on the fractional differential equations (FDEs). Fractional calculus is the generalizations of the ordinary differentiations and integrations to non-integer one. There are two kinds of fractional derivatives. The first kind is a nonlocal fractional derivative that is Caputo and Riemann–Liouville derivatives which applied successfully in different fields of science and engineering [7-13]. The second kind is the local fractional derivative, i.e., modified Riemann–Liouville (mRL), Kolwankar–Gangal (K–G), Cresson's, and Chen's fractal derivatives [14-23]. The purpose of using FDEs which is a powerful tool for the description of memory and hereditary properties of various materials and processes.

Nowadays studying FDEs become a hot topic of research owing to the development of the fractal theory and mathematical software package. There are many effective methods to treat numerical and analytical solutions of FDEs such as Adomian decomposition method, homotopy perturbation method, variational iteration method, finite difference method, Laplace transform, Fourier transform, generalized differential transform, Lie symmetry group, the fractional sub-equation method, the (*G′/G*)-expansion method, and the first integral method [24-41].

In the present paper, we introduce the fractional version of Lane-Emden type equations which play a vital role in physics and astrophysics. The main goal of this paper is to study several examples of Lane-Emden type by using series expansion method in the sense of mRL derivatives.

The structure of this paper is as follows: In section 2, we introduce some basic definitions of the fractional calculus theory. Several examples are solved in section 3. Finally, conclusions and discussions are given in section 4.

2. **Basic Definitions**

The mRL derivative definition [14-18]. Suppose that $f: R \to R, \quad x \to f(x)$ where $f(x)$ a continuous function, and let $h$ is a discretization span constant, the limit form of the mRL derivative

$$f^{(\alpha)}(x) = \lim_{h \downarrow 0} \frac{\Delta^{\alpha}[f(x) - f(0)]}{h^{\alpha}}, \qquad 0 < \alpha < 1,$$

where $\Delta^\alpha f(x) = \sum_{k=0}^{\infty} (-1)^k \frac{\Gamma(\alpha+1)}{\Gamma(k+1)\Gamma(\alpha-k+1)} f[x+(\alpha-k)h]$. Which is similar to the standard of derivatives, and the $\alpha$-order derivative of a constant is zero. The integral form of the mRL derivative is written as

$$D_x^\alpha f(x) = \begin{cases} \dfrac{1}{\Gamma(-\alpha)} \int_0^x (x-\xi)^{-\alpha-1} [f(\xi)-f(0)] d\xi, & \alpha < 0 \\ \dfrac{1}{\Gamma(1-\alpha)} \dfrac{d}{dx} \int_0^x (x-\xi)^{-\alpha} [f(\xi)-f(0)] d\xi, & 0 < \alpha < 1 \\ \dfrac{1}{\Gamma(n-\alpha)} \dfrac{d^n}{dx^n} \int_0^x (x-\xi)^{n-\alpha-1} [f(\xi)-f(0)] d\xi, & n \leq \alpha < n+1, \ n \geq 1. \end{cases} \quad (1)$$

Some properties of the mRL derivative are

$$D_x^\alpha x^\gamma = \frac{\Gamma(\gamma+1)}{\Gamma(\gamma+1-\alpha)} x^{\gamma-\alpha}, \qquad \gamma > 0, \tag{2}$$

$$D_x^\alpha (c f(x)) = c D_x^\alpha f(x), \tag{3}$$

$$D_x^\alpha [f(x)g(x)] = g(x) D_x^\alpha f(x) + f(x) D_x^\alpha g(x), \tag{4}$$

$$D_x^\alpha f[g(x)] = f_g'[g(x)] D_x^\alpha g(x), \tag{5}$$

$$D_x^\alpha f[g(x)] = D_g^\alpha f[g(x)] (g_x')^\alpha, \tag{6}$$

where $c$ is constant. The last five equations are direct results of $D_x^\alpha f(x) \cong \Gamma(\alpha+1) D_x f(x)$.

Based on the definition of the fractal index ($\sigma_x$) which is usually determined in terms of gamma functions [42- 44], Equations (4)-(6) could be modified to the following

$$D_x^\alpha [f(x)g(x)] = \sigma_x \{g(x) D_x^\alpha f(x) + f(x) D_x^\alpha g(x)\}, \tag{7}$$

$$D_x^\alpha f[g(x)] = \sigma_x f_g'[g(x)] D_x^\alpha g(x), \tag{8}$$

$$D_x^\alpha f[g(x)] = \sigma_x D_g^\alpha f[g(x)] (g_x')^\alpha. \tag{9}$$

In the next section, we discuss seven examples based on the fractal index and series expansion method.

**2. Series solution of class of Lane-Emden type equations**

**Example 1.** Consider the following fractional differential linear, homogeneous Lane-Emden equation

$$D_x^{\alpha\alpha} y + \frac{2}{x^\alpha} D_x^\alpha y - 2(4x^{2\alpha} + 3) y = 0, \qquad (10)$$

with initial conditions $y(0) = 1$, $D_x^\alpha y(0) = 0$, where $D_x^{\alpha\alpha} = D_x^\alpha D_x^\alpha$, $D_x^\alpha$ is the mRL-derivative, $y = y(x)$ is an unknown function and $0 < \alpha < 1$.

We assume the transform $X = x^\alpha$ and the solution can be expressed in a series in the form

$$y(X) = \sum_{m=0}^\infty A_m X^m, \qquad (11)$$

since $y(0) = 1$, $D_x^\alpha y(0) = 0$, this gives $A_0 = 1, A_1 = 0$. So $y = 1 + \sum_{m=2}^\infty A_m X^m$. Based on the fractal index we can calculate the $\alpha$-fractional derivative and the twice $\alpha$-fractional derivative with respect to $x$

$$D_x^\alpha y = \sum_{m=2}^\infty A_m \sigma_{1x} m X^{m-1} \frac{\Gamma(\alpha+1)}{\Gamma(\alpha+1-\alpha)} x^{\alpha-\alpha}$$

where $\sigma_{1x}$ is called the fractal index [38] given by

$$\sigma_{1x} = \frac{\Gamma(m\alpha+1)}{m\Gamma(\alpha+1)\Gamma(m\alpha+1-\alpha)},$$

so

$$D_x^\alpha y = \sum_{m=2}^\infty A_m X^{m-1} \frac{\Gamma(m\alpha+1)}{\Gamma(m\alpha+1-\alpha)}, \qquad (12)$$

$$D_x^{\alpha\alpha} u = \sum_{m=2}^\infty A_m \sigma_{2x} (m-1) X^{m-2} \frac{\Gamma(m\alpha+1)}{\Gamma(m\alpha+1-\alpha)} \frac{\Gamma(\alpha+1)}{\Gamma(\alpha+1-\alpha)} x^{\alpha-\alpha} = \sum_{m=2}^\infty X^{m-2} \frac{A_m \Gamma(m\alpha+1)}{\Gamma(m\alpha+1-2\alpha)}, \qquad (13)$$

where

$$\sigma_{2x} = \frac{\Gamma((m-1)\alpha+1)}{(m-1)\Gamma(\alpha+1)\Gamma((m-1)\alpha+1-\alpha)}.$$

Substituting Equations (11), (12) and (13) into Equation (10), yields

$$\sum_{m=2}^\infty X^m \frac{A_m \Gamma(m\alpha+1)}{\Gamma(m\alpha+1-2\alpha)} + \sum_{m=2}^\infty A_m X^m \frac{2\Gamma(m\alpha+1)}{\Gamma(m\alpha+1-\alpha)} - 6X^2 - 4X^4 - \sum_{m=4}^\infty 6A_{m-2} X^m - \sum_{m=6}^\infty 4A_{m-4} X^m = 0.$$

Equating the coefficients of $X^2, X^3, X^4, X^5$ and $X^m$ we get

$$A_2 = \frac{6\Gamma(\alpha+1)}{\Gamma(2\alpha+1)[\Gamma(\alpha+1)+2]}, \qquad A_3 = 0,$$

$$A_4 = \frac{\Gamma(2\alpha+1)\Gamma(3\alpha+1)}{\Gamma(4\alpha+1)[\Gamma(3\alpha+1)+2\Gamma(2\alpha+1)]}\left[4+\frac{36\Gamma(\alpha+1)}{\Gamma(2\alpha+1)[\Gamma(\alpha+1)+2]}\right], \qquad A_5 = 0,$$
(14)

$$\frac{A_m\Gamma(m\alpha+1)}{\Gamma(m\alpha+1-2\alpha)} + \frac{2A_m\Gamma(m\alpha+1)}{\Gamma(m\alpha+1-\alpha)} - 6A_{m-2} - 4A_{m-4} = 0, \quad m = 6,7,8,\ldots. \tag{15}$$

Equation (15) is the recurrence relation of the series expansion, Equation (10). The series solution of Equation (10) is

$$y = 1 + \frac{6\Gamma(\alpha+1)x^{2\alpha}}{\Gamma(2\alpha+1)[\Gamma(\alpha+1)+2]} + \left[4+\frac{36\Gamma(\alpha+1)}{\Gamma(2\alpha+1)[\Gamma(\alpha+1)+2]}\right]\frac{\Gamma(2\alpha+1)\Gamma(3\alpha+1)x^{4\alpha}}{\Gamma(4\alpha+1)[\Gamma(3\alpha+1)+2\Gamma(2\alpha+1)]} + \ldots$$
(16)

When $\alpha = 1$, Equation (10) becomes

$$y'' + \frac{2}{x}y' - 2(4x^2+3)y = 0, \tag{17}$$

and has a series solution in the form

$$y = 1 + x^2 + \frac{x^4}{2!} + \frac{x^6}{3!} + \ldots = e^{x^2} \tag{18}$$

which are the same results obtained using Adomian decomposition method [45], homotopy perturbation method [46], modified Laplace decomposition method [47] and variational iteration method [48].

**Example 2.** Consider the following fractional differential linear, homogeneous Lane-Emden equation

$$D_x^{\alpha\alpha}y + \frac{2}{x^\alpha}D_x^\alpha y + y^n = 0, \quad x \geq 0, \quad n = 0,1, \tag{19}$$

with initial conditions $y(0) = 1$, $D_x^\alpha y(0) = 0$.

When $n = 0$, Substituting Equations (12) and (13) into Equation (19), yields

$$\sum_{m=2}^{\infty} X^m \frac{A_m\Gamma(m\alpha+1)}{\Gamma(m\alpha+1-2\alpha)} + \sum_{m=2}^{\infty} A_m X^m \frac{2\Gamma(m\alpha+1)}{\Gamma(m\alpha+1-\alpha)} + X^2 = 0.$$

Equating the coefficients of $X^2$ and $X^m$ we get

$$A_2 = -\frac{\Gamma(\alpha+1)}{\Gamma(2\alpha+1)[\Gamma(\alpha+1)+2]}, \quad A_m = 0, \quad m = 3,4,5,\dots. \tag{20}$$

Thus the solution has the form

$$y = 1 - \frac{\Gamma(\alpha+1)x^{2\alpha}}{\Gamma(2\alpha+1)[\Gamma(\alpha+1)+2]}. \tag{21}$$

When $\alpha = 1$, solution of Equation (21) becomes

$$y = 1 - \frac{x^2}{6},$$

which is the exact solution of the Lane-Emden equation [45-52] at $n = 0$.

When $n = 1$, substituting Equations (11), (12) and (13) into Equation (19), yields

$$\sum_{m=2}^{\infty} X^m \frac{A_m \Gamma(m\alpha+1)}{\Gamma(m\alpha+1-2\alpha)} + \sum_{m=2}^{\infty} A_m X^m \frac{2\Gamma(m\alpha+1)}{\Gamma(m\alpha+1-\alpha)} + X^2 + \sum_{m=4}^{\infty} A_{m-2} X^m = 0,$$

Equating the coefficients of $X^2, X^3$ and $X^m$ we get

$$A_2 = -\frac{\Gamma(\alpha+1)}{\Gamma(2\alpha+1)[\Gamma(\alpha+1)+2]}, \qquad A_3 = 0,$$

$$\frac{A_m \Gamma(m\alpha+1)}{\Gamma(m\alpha+1-2\alpha)} + \frac{2A_m \Gamma(m\alpha+1)}{\Gamma(m\alpha+1-\alpha)} + A_{m-2} = 0, \quad m = 4,5,6,\dots, \tag{22}$$

The series solution takes the form

$$y = 1 - \frac{\Gamma(\alpha+1)x^{2\alpha}}{\Gamma(2\alpha+1)[\Gamma(\alpha+1)+2]} + \frac{\Gamma(\alpha+1)\Gamma(3\alpha+1)x^{4\alpha}}{\Gamma(4\alpha+1)[\Gamma(3\alpha+1)+2\Gamma(2\alpha+1)][\Gamma(\alpha+1)+2]} + \dots. \tag{23}$$

When $\alpha = 1$, Equation (23) will be

$$y = 1 - \frac{x^2}{3!} + \frac{x^4}{5!} - \frac{x^6}{7!} + \dots = \frac{\sin x}{x}.$$

which is the well-known solution of the Lane-Emden Equation [45-52] at $n = 1$.

**Example 3.** Consider the following fractional differential linear, homogeneous Lane-Emden equation

$$D_x^{\alpha\alpha} y + \frac{2}{x^\alpha} D_x^\alpha y + y = 6 + 12x^\alpha + x^{2\alpha} + x^{3\alpha}, \quad 0 < x \le 2, \tag{24}$$

with initial conditions $y(0) = 0, \; D_x^\alpha y(0) = 0.$

From the initial conditions we have $A_0 = 0, A_1 = 0$. So $y = \sum_{m=2}^{\infty} A_m X^m$. Substituting Equations (12) and (13) into equation (24), yields

$$\sum_{m=2}^{\infty} X^m \frac{A_m \Gamma(m\alpha+1)}{\Gamma(m\alpha+1-2\alpha)} + \sum_{m=2}^{\infty} A_m X^m \frac{2\Gamma(m\alpha+1)}{\Gamma(m\alpha+1-\alpha)} + \sum_{m=4}^{\infty} A_{m-2} X^m - 6X^2 - 12X^3 - X^4 - X^5 = 0. \quad (25)$$

Equating the coefficient $X^m$ we get the general expression of the series coefficients as

$$\frac{\Gamma(m\alpha+1)[\Gamma((m-1)\alpha+1) + 2\Gamma((m-2)\alpha+1)]}{\Gamma((m-1)\alpha+1)\Gamma((m-2)\alpha+1)} A_m + A_{m-2} = 0, \quad m = 6, 7, ....$$

The series solution of Equation (24) could be given by

$$y = \frac{6\Gamma(\alpha+1)x^{2\alpha}}{\Gamma(2\alpha+1)[\Gamma(\alpha+1)+2]} + \frac{12\Gamma(\alpha+1)\Gamma(2\alpha+1)x^{3\alpha}}{\Gamma(3\alpha+1)[\Gamma(2\alpha+1)+2\Gamma(\alpha+1)]}$$
$$+ \frac{\Gamma(3\alpha+1)\Gamma(2\alpha+1)\{\Gamma(2\alpha+1)[\Gamma(\alpha+1)+2] - 6\Gamma(\alpha+1)\}x^{4\alpha}}{\Gamma(2\alpha+1)[\Gamma(\alpha+1)+2]\Gamma(4\alpha+1)[\Gamma(3\alpha+1)+2\Gamma(2\alpha+1)]} + \ldots . \quad (26)$$

When $\alpha = 1$, the solution of Equation (26) will be $y = x^2 + x^3$, [45-48].

**Example 4.** Consider the following fractional differential nonlinear Lane-Emden equation

$$D_x^{\alpha\alpha} y + \frac{2}{x^\alpha} D_x^\alpha y - y^n = 0, \quad x \geq 0, \quad (27)$$

with initial conditions $y(0) = 1$, $D_x^\alpha y(0) = 0$.

Suppose that

$$G(X) = \sum_{m=0}^{\infty} Q_m X^m = Q_0 + Q_1 X + Q_2 X^2 + Q_3 X^3 + Q_4 X^4 + Q_5 X^5 + \ldots \quad (28)$$

and

$$y^n(X) = G(X) \quad (29)$$

To obtain the fractional derivative of $y^n$, we apply the fractional derivative of the product of two functions, since the fractional derivative of $y^2$ will be considered as $y$ times $y$. Similarly $y^3$ will be considered as $y$ times $y^2$ and so on. Taking the fractional derivative for the both sides of Equation (29), we have

$$n y^{n-1} D_x^\alpha y = D_x^\alpha G \quad \text{or} \quad n y^n D_x^\alpha y = y D_x^\alpha G, \quad (30)$$

Differentiating both sides of Equation (30) $k$ times $\alpha$-derivatives we have

$$n\sum_{j=0}^{k}\binom{k}{j}\underbrace{D_x^{\alpha}...D_x^{\alpha}}_{(j+1)times}u.\underbrace{D_x^{\alpha}...D_x^{\alpha}}_{(k-j)times}G = \sum_{j=0}^{k}\binom{k}{j}\underbrace{D_x^{\alpha}...D_x^{\alpha}}_{j+1 times}G.\underbrace{D_x^{\alpha}...D_x^{\alpha}}_{k-j times}u, \qquad (31)$$

at $x = 0$, we have

$$n\sum_{j=0}^{k}\binom{k}{j}\underbrace{D_x^{\alpha}...D_x^{\alpha}}_{(j+1)times}u(0).\underbrace{D_x^{\alpha}...D_x^{\alpha}}_{(k-j)times}G(0) = \sum_{j=0}^{k}\binom{k}{j}\underbrace{D_x^{\alpha}...D_x^{\alpha}}_{j+1 times}G(0).\underbrace{D_x^{\alpha}...D_x^{\alpha}}_{k-j times}u(0) \qquad (32)$$

since

$$\underbrace{D_x^{\alpha}...D_x^{\alpha}}_{(j+1)times}u(0) = A_{j+1}\Gamma((j+1)\alpha+1), \qquad \underbrace{D_x^{\alpha}...D_x^{\alpha}}_{(k-j)times}G(0) = Q_{k-j}\Gamma(\alpha(k-j)+1),$$

$$\underbrace{D_x^{\alpha}...D_x^{\alpha}}_{(j+1)times}G(0) = Q_{j+1}\Gamma((j+1)\alpha+1), \qquad \underbrace{D_x^{\alpha}...D_x^{\alpha}}_{(k-j)times}u(0) = A_{k-j}\Gamma(\alpha(k-j)),$$

making some algebraic manipulations we get

$$\Gamma(m\alpha+1)A_0Q_m = n\sum_{i=1}^{m}\frac{(m-1)!\Gamma(\alpha(m-i)+1)\Gamma(i\alpha+1)}{(i-1)!(m-i)!}A_iQ_{m-i}$$
$$- \sum_{i=1}^{m}\frac{(m-1)!\Gamma((m-i)\alpha+1)\Gamma(i\alpha+1)}{(m-1-i)!i!}A_iQ_{m-i}, \qquad (33)$$

the recurrence relation could be written as

$$Q_m = \frac{1}{\Gamma(m\alpha+1)A_0}\sum_{i=1}^{m}\frac{(m-1)!\Gamma(\alpha(m-i)+1)\Gamma(i\alpha+1)}{i!(m-i)!}[in-m+i]A_iQ_{m-i}, \quad \forall\ m \geq 1, \qquad (34)$$

with $A_0 = 1$, $A_1 = 0$, we have the two coefficients

$$Q_0 = A_0^n = 1, \qquad Q_1 = \frac{\Gamma(\alpha+1)n}{\Gamma(\alpha+1)A_0}A_1Q_0 = 0.$$

Substituting Equations (3), (4) and (20) into Equation (18), yields

$$\sum_{m=2}^{\infty}X^m\frac{A_m\Gamma(m\alpha+1)}{\Gamma(m\alpha+1-2\alpha)} + \sum_{m=2}^{\infty}X^m\frac{2A_m\Gamma(m\alpha+1)}{\Gamma(m\alpha+1-\alpha)} - X^2 - \sum_{m=4}^{\infty}Q_{m-2}X^m = 0. \qquad (35)$$

Equating the coefficients of $X^2, X^3$ and $X^m$ we get

$$A_2 = \frac{\Gamma(\alpha+1)}{\Gamma(2\alpha+1)[\Gamma(\alpha+1)+2]}, \qquad A_3 = 0$$

$$\frac{A_{m+2}\Gamma((m+2)\alpha+1)}{\Gamma(m\alpha+1)} + \frac{2A_{m+2}\Gamma((m+2)\alpha+1)}{\Gamma((m+1)\alpha+1)} - Q_m = 0, \quad , m = 2,3,4,.... \qquad (36)$$

The series solution of equation (27)

$$y = 1 + \frac{\Gamma(\alpha+1) \, x^{2\alpha}}{\Gamma(2\alpha+1)[\Gamma(\alpha+1)+2]} + \frac{n\Gamma(\alpha+1)\Gamma(2\alpha+1)\Gamma(3\alpha+1) \, x^{4\alpha}}{\Gamma(4\alpha+1)\Gamma(2\alpha+1)[\Gamma(3\alpha+1)+2\Gamma(2\alpha+1)][\Gamma(\alpha+1)+2]} + \ldots \quad (37)$$

From the above procedure, we can get the solution of the fractional differential nonlinear Lane-Emden equation

$$D_x^{\alpha\alpha} y + \frac{2}{x^\alpha} D_x^\alpha y + y^n = 0, \quad x \geq 0, \quad (38)$$

with initial conditions $y(0) = 1$, $D_x^\alpha y(0) = 0$, in the following form

$$y = 1 - \frac{\Gamma(\alpha+1) \, x^{2\alpha}}{\Gamma(2\alpha+1)[\Gamma(\alpha+1)+2]} + \frac{n\Gamma(\alpha+1)\Gamma(2\alpha+1)\Gamma(3\alpha+1) \, x^{4\alpha}}{\Gamma(4\alpha+1)\Gamma(2\alpha+1)[\Gamma(3\alpha+1)+2\Gamma(2\alpha+1)][\Gamma(\alpha+1)+2]} - \ldots \quad (39)$$

When $\alpha = 1$, solutions (37) and (38) will be

$$y_- = 1 + \frac{x^2}{6} + \frac{n\,x^4}{120} + \ldots, \quad (40)$$

$$y_+ = 1 - \frac{x^2}{6} + \frac{n\,x^4}{120} - \ldots \quad (41)$$

which are similar to [6].

**Example 5.** Consider the following fractional differential nonlinear Lane-Emden equation

$$D_x^{\alpha\alpha} y + \frac{2}{x^\alpha} D_x^\alpha y + 4 x^{2\alpha} y^2 = 0, \quad x \geq 0, \quad (42)$$

with initial conditions $y(0) = 1$, $D_x^\alpha y(0) = 0$.

From the initial conditions we have $A_0 = 1, A_1 = 0$. So $y = 1 + \sum_{m=2}^{\infty} A_m X^m$. Substituting Equations (12) and (13) into Equation (32) with the recurrence relation (24) at $n = 2$, yields

$$\sum_{m=2}^{\infty} X^m \frac{A_m \Gamma(m\alpha+1)}{\Gamma(m\alpha+1-2\alpha)} + \sum_{m=2}^{\infty} X^m \frac{2A_m \Gamma(m\alpha+1)}{\Gamma(m\alpha+1-\alpha)} + 4X^2 + \sum_{m=4}^{\infty} 4Q_{m-2} X^m = 0. \quad (43)$$

Equating the coefficients of $X^2, X^3$ and $X^m$ we get

$$A_2 = \frac{4\Gamma(\alpha+1)}{\Gamma(2\alpha+1)[\Gamma(\alpha+1)+2]}, \qquad A_3 = 0$$

$$\frac{A_{m+2}\Gamma((m+2)\alpha+1)}{\Gamma(m\alpha+1)} + \frac{2A_{m+2}\Gamma((m+2)\alpha+1)}{\Gamma((m+1)\alpha+1)} + 4Q_m = 0, \quad m = 2,3,4,\ldots.$$

From Equation (34), since $A_0 = 1, A_1 = 0$ we have $Q_0 = 1, Q_1 = 0, Q_2 = \dfrac{8\Gamma(\alpha+1)}{\Gamma(2\alpha+1)[\Gamma(\alpha+1)+2]}$, this gives $A_4 = -\dfrac{32\Gamma(\alpha+1)\Gamma(3\alpha+1)}{\Gamma(4\alpha+1)[\Gamma(3\alpha+1)+2\Gamma(2\alpha+1)][\Gamma(\alpha+1)+2]}$ and so on. The series solution of

Equation (42) could be given by

$$y = 1 + \frac{4\Gamma(\alpha+1)\, x^{2\alpha}}{\Gamma(2\alpha+1)[\Gamma(\alpha+1)+2]} - \frac{32\Gamma(\alpha+1)\Gamma(3\alpha+1)\, x^{4\alpha}}{\Gamma(4\alpha+1)[\Gamma(3\alpha+1)+2\Gamma(2\alpha+1)][\Gamma(\alpha+1)+2]} + \dots \quad (44)$$

When $\alpha = 1$, solution of Equation (44) will be given by $y = 1 + \dfrac{2}{3}x^2 - \dfrac{4}{15}x^4 + \dots$ .

**Example 6.** Consider the following fractional differential nonlinear Lane-Emden equation

$$D_x^{\alpha\alpha} y + \frac{2}{x^\alpha} D_x^\alpha y + y^3 = 6 + x^{6\alpha}, \quad x \geq 0, \quad (45)$$

with initial conditions $y(0) = 0$, $D_x^\alpha y(0) = 0$.

From the initial conditions we have $A_0 = 0, A_1 = 0$. So $y = \sum\limits_{m=2}^{\infty} A_m X^m$. Substituting

Equations (12) and (13) into Equation (45) with the recurrence relation (34) at $n = 3$, yields

$$\sum_{m=2}^{\infty} X^m \frac{A_m \Gamma(m\alpha+1)}{\Gamma(m\alpha+1-2\alpha)} + \sum_{m=2}^{\infty} X^m \frac{2 A_m \Gamma(m\alpha+1)}{\Gamma(m\alpha+1-\alpha)} + \sum_{m=4}^{\infty} Q_{m-2} X^m = 6X^2 + X^8. \quad (46)$$

Equating the coefficients of $X^2, X^3, X^4, X^5, X^6, X^7, X^8$ and $X^m$ we get

$$A_2 = \frac{6\Gamma(\alpha+1)}{\Gamma(2\alpha+1)[\Gamma(\alpha+1)+2]}, \quad A_3 = 0$$

$$\frac{A_{m+2}\Gamma((m+2)\alpha+1)}{\Gamma(m\alpha+1)} + \frac{2A_{m+2}\Gamma((m+2)\alpha+1)}{\Gamma((m+1)\alpha+1)} + Q_m = 0, \quad m = 2,3,4,5,7,8,\dots,$$

$$\frac{A_8\Gamma((m+2)\alpha+1)}{\Gamma(m\alpha+1)} + \frac{2A_8\Gamma((m+2)\alpha+1)}{\Gamma((m+1)\alpha+1)} + Q_6 = 1.$$

From Equation (34), since $A_0 = 0, A_1 = 0$ we have $Q_0 = 0, Q_1 = 0, Q_2 = 0$, this gives $A_4 = 0$. Also $Q_3 = 0$ gives $A_5 = 0$, $Q_4 = 0$ gives $A_6 = 0$, $Q_5 = 0$ gives $A_7 = 0$ and $Q_6 = 0$ gives

$$A_8 = \frac{\Gamma(6\alpha+1)\Gamma(7\alpha+1)}{\Gamma(8\alpha+1)[\Gamma(7\alpha+1)+2\Gamma(6\alpha+1)]},$$ and so on. The series solution of Equation (45)

$$y = \frac{6\Gamma(\alpha+1)\, x^{2\alpha}}{\Gamma(2\alpha+1)[\Gamma(\alpha+1)+2]} + \frac{\Gamma(6\alpha+1)\Gamma(7\alpha+1)\, x^{8\alpha}}{\Gamma(8\alpha+1)[\Gamma(7\alpha+1)+2\Gamma(6\alpha+1)]}. \quad (47)$$

When $\alpha = 1$, solution (47) will be $y = x^2 + \dfrac{1}{72}x^8$.

**Example 7.** Consider the fractional Lane-Emden equation of the isothermal gas sphere

$$D_x^{\alpha\alpha} y + \frac{2}{x^\alpha} D_x^\alpha y - e^{-y} = 0, \quad x \geq 0, \tag{48}$$

with initial conditions $y(0) = 0$, $D_x^\alpha y(0) = 0$. Equation (48) describing gas spheres, isothermal gas spheres embedded in a pressurized medium at the maximum possible mass allowing for hydrostatic equilibrium [6].

From the initial conditions we have $A_0 = 0, A_1 = 0$. So $y = \sum_{m=2}^{\infty} A_m X^m$. By putting

$$e^{-y} = S(X) = \sum_{m=0}^{\infty} R_m X^m, \tag{49}$$

we have $e^{-y(0)} = S(0) = 1$, or $R_0 = 1$. Doing the same procedures as equation (20) we have $D_x^\alpha e^{-y} = D_x^\alpha S$, gives

$$-G D_x^\alpha y = D_x^\alpha S, \tag{50}$$

Differentiating both sides of Equation (50) $k$ times $\alpha$-derivatives we have

$$\sum_{j=0}^{k} \binom{k}{j} \underbrace{D_x^\alpha \ldots D_x^\alpha}_{(j+1)\text{times}} y \cdot \underbrace{D_x^\alpha \ldots D_x^\alpha}_{(k-j)\text{times}} S = -\underbrace{D_x^\alpha \ldots D_x^\alpha}_{j+1\text{times}} S,$$

at $x = 0$, we have

$$\sum_{j=0}^{k} \binom{k}{j} \underbrace{D_x^\alpha \ldots D_x^\alpha}_{(j+1)\text{times}} y(0) \cdot \underbrace{D_x^\alpha \ldots D_x^\alpha}_{(k-j)\text{times}} S(0) = -\underbrace{D_x^\alpha \ldots D_x^\alpha}_{j+1\text{times}} S(0),$$

we have

$$R_l = -\frac{(l-1)!}{\Gamma(l\alpha+1)} \sum_{i=1}^{l} \frac{A_i \Gamma(i\alpha+1) R_{l-i} \Gamma((l-i)\alpha+1)}{(i-1)!(l-i)!}. \tag{51}$$

Substituting Equations (12), (13) and (51) into Equation (35), yields

$$\sum_{m=0}^{\infty} X^{m+2} \frac{A_{m+2}\Gamma((m+2)\alpha+1)}{\Gamma(m\alpha+1)} + \sum_{m=0}^{\infty} X^{m+2} \frac{2A_{m+2}\Gamma((m+2)\alpha+1)}{\Gamma((m+1)\alpha+1)} - X^2 - \sum_{m=1}^{\infty} R_m X^{m+2} = 0. \tag{52}$$

Equating the coefficients of $X^2$ and $X^{m+2}$ we get

$$A_2 = \frac{\Gamma(\alpha+1)}{\Gamma(2\alpha+1)[\Gamma(\alpha+1)+2]}, \tag{53}$$

$$\frac{A_{m+2}\Gamma((m+2)\alpha+1)}{\Gamma(m\alpha+1)} + \frac{2A_{m+2}\Gamma((m+2)\alpha+1)}{\Gamma((m+1)\alpha+1)} - R_m = 0, \tag{54}$$

The series solution of Equation (48) is

$$y = \frac{\Gamma(\alpha+1)\, x^{2\alpha}}{\Gamma(2\alpha+1)[\Gamma(\alpha+1)+2]} - \frac{\Gamma(\alpha+1)\Gamma(3\alpha+1)\, x^{4\alpha}}{\Gamma(4\alpha+1)[\Gamma(3\alpha+1)+2\Gamma(2\alpha+1)][\Gamma(\alpha+1)+2]}$$
$$+ \frac{\Gamma(\alpha+1)\Gamma(5\alpha+1)\, x^{6\alpha}}{\Gamma(6\alpha+1)[\Gamma(5\alpha+1)+2\Gamma(4\alpha+1)][\Gamma(\alpha+1)+2]} \left[ \frac{3\Gamma(\alpha+1)}{\Gamma(\alpha+1)+2} + \frac{\Gamma(3\alpha+1)}{\Gamma(3\alpha+1)+2\Gamma(2\alpha+1)} \right] + \ldots \tag{55}$$

When $\alpha = 1$, solution of Equation (55) goes to the will known solution of the isothermal gas sphere equation [6]

$$y = \frac{1}{6}x^2 - \frac{1}{120}x^4 + \frac{1}{1890}x^6 + \ldots . \tag{56}$$

## 4. Conclusions

Based on the series method, approximate solutions can be easily obtained for a class of fractional Lane–Emden type equations with the sense of the modified Riemann–Liouville derivative. The fractional Lane–Emden type equation which can model many phenomena in mathematical physics and astrophysics. The desired algorithm is reliable and easy to use. Its validity is verified by various examples. When $\alpha = 1$, our solutions recover the corresponding results given by Adomian decomposition method, homotopy perturbation method, modified Laplace decomposition method and variational iteration method. The results show that the series method is very effective and convenient for solving the fractional Lane-Emden type equations. The series method can be applied to other linear and nonlinear FDEs.

[49] M. Hussain and M. Khan, *Applied Mathematical Sciences* **4** (2010) 1769.

[50] S. A. Yousefi, *Appl. Math. Comput.* **181**(2006) 1417.

[51] J. I. Ramos, *Appl. Math. Comput.* **161** (2005) 525.15